\documentclass[10pt]{amsart}

\usepackage{amsfonts}
\usepackage{latexsym}
\usepackage{color}
 \RequirePackage{amsbsy} 
 \RequirePackage{amsopn} 
 \RequirePackage{amsmath}
  \RequirePackage{amssymb}
  \usepackage[mathscr]{eucal}
\definecolor{darkmagenta}{rgb}{0.5, 0, 0.5}
\definecolor{darkblue}{rgb}{0.1, 0.1, 0.7}
\definecolor{darkgreen}{rgb}{0.1, 0.35, 0.1}

 \newtheorem{thee}{Theorem}
 
 \newtheorem{leem}[thee]{Lemma}
 
 \newtheorem{prro}[thee]{Proposition}
 \newtheorem{deef}[thee]{Definition}

 \newcommand{\balf}
 {\renewcommand{\theenumi}{(\alph{enumi})}
 \renewcommand{\labelenumi}{\theenumi}
                      \begin{enumerate}}
\newcommand{\ealf}   {\end{enumerate}
                      \renewcommand{\theenumi}{\arabic{enumi}}
                      \renewcommand{\labelenumi}{\theenumi.}}
\newcommand{\bara}   {\renewcommand{\theenumi}{(\arabic{enumi})}
                      \renewcommand{\labelenumi}{\theenumi}
                      \begin{enumerate} }
\newcommand{\eara}   {\end{enumerate}
                      \renewcommand{\theenumi}{\arabic{enumi}}
                      \renewcommand{\labelenumi}{\theenumi.}}

 \newcommand{\brom}   {\renewcommand{\theenumi}{(\roman{enumi})}
                      \renewcommand{\labelenumi}{\theenumi}
                      \begin{enumerate} }
\newcommand{\erom}   {\end{enumerate}
                      \renewcommand{\theenumi}{\arabic{enumi}}
                      \renewcommand{\labelenumi}{\theenumi.}}

	   \newcommand{\Spec}{\mbox{\rm Spec}}

  \DeclareMathOperator{\C}     {\boldsymbol C}%
\DeclareMathOperator{\Z}     {\boldsymbol Z}%
  \DeclareMathOperator{\calU} {\mathcal U}%
  \RequirePackage{amssymb} 
  \RequirePackage{amsmath} 

\begin{document}

 \title[patch topology and ultrafilter topology]{The patch topology 
and the ultrafilter topology on the prime spectrum of a commutative 
ring}

   \author[ ]{Marco Fontana$^{1}$ \; and\ \; K. Alan Loper$^{2}$}

 \maketitle

  \dedicatory{\rm \begin{center} $^{1}$ Dipartimento di Matematica \\
  Universit\`a degli Studi Roma Tre \\
  Largo San Leonardo Murialdo, 1 \\ 
  00146 Roma, Italy\\ {\texttt{fontana@mat.uniroma3.it}} \\
  \texttt{              } \\ \rm $^{2}$ Department of Mathematics \\ 
  Ohio State
  University-Newark\\ Newark, Ohio 43055 USA \\
  {\texttt{lopera@math.ohio-state.edu}} \rm \end{center}}

  \bigskip
  \bigskip
  \bigskip

   \begin{abstract}
Let $R$ be a commutative ring and let $\Spec(R)$ denote the collection 
of prime ideals of $R$.  We define a topology on $\Spec(R)$ by using 
ultrafilters and demonstrate that this topology is identical to the 
well  known {\it patch } or {\it constructible} topology.  The proof 
is accomplished by use of a von Neumann regular ring canonically associated with $R$.
 \end{abstract}

  \bigskip





Let $R$ be a commutative ring and let $\Spec(R)$ denote the 
collection 
of prime ideals of $R$.  On $\Spec(R)$ we can define a topology 
known as  {\it Zariski's topology}:
the collection of all sets $V(I):= \{P \in \Spec(R) \, | \, I 
\subseteq P \}$ where $I$ is an ideal of $R$ constitutes the closed sets in 
this topology.

Zariski's topology has several attractive properties related to the 
geometric aspects of the study of the set of prime ideals 
\cite[Chapter I]{Eisenbud}.  For example, $\Spec(R)$ is 
always quasicompact (that is, every open covering has a finite 
refinement).  On the other hand, this topology is very coarse.  For 
example, 
$\Spec(R)$ is almost never 
Hausdorff (that is, two distinct points have nonintersecting 
neighborhoods).  

Many authors have considered a finer 
topology, known as the {\it patch topology}  \cite{Hochster} and as 
the 
{\it constructible topology} (\cite [pages 337-339]{EGA} or  
\cite[Chapter 3,  Exercises 27, 28 and 30]{AM}), which can be defined 
starting from Zariski's topology.  

Consider two collections of subsets of $\Spec(R)$.

\begin{enumerate}
\item The sets $V(I)$ defined above for $I$ an ideal of  $R$.
\item The sets $D(a) := \Spec(R) \setminus V(a)$ where $a\in R$   
(where, as usual,  $V(a)$ denotes the set $V(aR)$). 

\end{enumerate}

The patch topology is then the smallest topology in which both of the 
above classes consist of closed sets.  The patch 
topology is a refinement of Zariski's topology which is always 
Hausdorff.

It is easy and natural to define another topology on $\Spec(R)$ by 
introducing the notion of an ultrafilter.

We start with definitions, some notation and some preliminary results.

Given an infinite set $S$  \it an ultrafilter $\calU$ on $S$ \rm 
 is a collection of 
subsets of $S$ such that

\begin{enumerate}
    
\item If $A \in \calU$ and $A \subseteq B \subseteq S$ then $B \in 
\calU$.
\item If $A,B \in \calU$ then $A \cap B \in \calU$.
\item If $A \cup B \in \calU$ and $A \cap B$ is empty then exactly 
one 
of $A $ and $B$ lies in $\calU$.

\end{enumerate}

Note that (1) implies that any nonempty ultrafilter on a set $S$ contains the set $S$ itself.  
It then follows from (3) that the empty set can never be a member of an ultrafilter.  Hence 
an ultrafilter on a set $S$ is always a proper subset of the power set of $S$.

 A principal ultrafilter on $S$ can be defined as follows. If $d 
\in S$ then  \it  the principal ultrafilter of $d$ on $S$ \rm   
is
the collection of all subsets of $S$ which contain $d$.  Zorn's Lemma 
can be used to prove that nonprincipal ultrafilters exist on any 
infinite  set.

Let $C$ be a subset of $\Spec(R)$, and 
let $\calU$ be an ultrafilter on the set $C$.   Set   $P_{\calU} 
:= \{a \in R  \mid V(a) \cap C \in \calU \}$.  By an argument similar 
to that used in \cite[Lemma 2.4]{CLT} it can be easily shown that 
$P_{\calU}$ is a prime ideal of $R$.     We 
call $P_{\calU}$ an \it ultrafilter limit point of \rm $C$. This 
notion of 
ultrafilter limits of collections of prime ideals has been used to 
great effect in 
several recent papers  \cite{CLT}, \cite{Loper1}, and  
\cite{Loper2}.  If $\calU$ is a principal ultrafilter then there is a 
prime $P \in C$ such that $\calU$ consists of all subsets of $C$ which 
contain $P$.  It is clear then that $P_{\calU} = P \in C$.  On the 
other hand, if $\calU$ is nonprincipal, then it is not at all 
clear that $P_{\calU}$ should lie in $C$.   That 
motivates our definition.

\begin{deef} \rm  Let $R$ and $C$ be as above.  We say that $C$ is 
\it ultrafilter closed \rm  if it contains all of its ultrafilter 
limit 
points. 
\end{deef}

It is not hard to see that the ultrafilter closed subsets 
of $\Spec(R)$ define a topology on the set $\Spec(R)$, called \it the 
ultrafilter topology on \rm $\Spec(R)$.   In fact:

\begin{itemize}

\item Suppose that $C_{1}, C_{2}, \ldots, C_{n}$ are ultrafilter 
closed subsets of $\Spec(R)$.  Let $C  := C_{1} \cup C_{2} \cup 
\ldots \cup C_{n}$.   
Let $\calU$ be an ultrafilter on $C$.  Then $\calU$ defines an 
ultrafilter limit prime $P_{\calU}$.  We want to show that 
$P_{\calU} \in C$.  Note that at least one    of the sets 
$C_{i}$ lies in $\calU$.  Without loss of  generality, suppose 
that 
$C_{1} \in \calU$.  The collection of sets 
$\calU_{1}   :=   \{ C_{1} \cap B \, | \, B \in \calU \}$ is 
then an 
ultrafilter on $C_{1}$ and the ultrafilter prime it defines is 
indentical to $P_{\calU}$.  In particular,

\begin{enumerate}
    
\item Let $d \in P_{\calU_{1}}$.  We know that 
$V(d) \cap C_{1} \in \calU_{1}$.  Since $C_{1} \in \calU$ then  
every set in $\calU_{1}$ is also in $\calU$.  Hence 
$V(d) \cap C_{1} \in \calU$.  Then note that 
$V(d) \cap C_{1} \subseteq V(d) \cap C$ and so 
$V(d) \cap C \in \calU$.  Hence, $d \in P_{\calU}$.

\item Let $d \in P_{\calU}$.  Then $V(d) \cap C \in \calU$.  The 
definition of $\calU_{1}$ then implies that 
$V(d) \cap C_{1} = (V(d) \cap C) \cap C_{1} \in \calU_{1}$.  Hence, 
$d \in P_{\calU_{1}}$.

\end{enumerate}

  Hence, $P_{\calU} \in C$ since 
it is in the ultrafilter closure of $C_{1}$ which is an ultrafilter  closed set. 

\item Suppose that  $ \{C_{\lambda} \, | \, \lambda \in \Lambda\}$  
is a collection of ultrafilter closed subsets of $\Spec(R)$.  Let
$C  :=   \cap_{\lambda \in \Lambda} C_{\lambda}$.  Let 
$\calU$ be an 
ultrafilter on $C$   and $P_{\calU}$ the  ultrafilter limit 
prime  associated to $\calU$. 
 We want to show that $P_{\calU}$ lies in $C$.  For 
each $\lambda \in \Lambda$ the collection 
$\calU_{\lambda}  :=  \{B \subseteq C_{\lambda} \mid B  
\cap C \in \calU \}$ is an ultrafilter on $C_{\lambda}$.  Moreover, it 
defines the 
same limit prime $P_{\calU}$ as the ultrafilter $\calU$ defines on 
$C$ (using an argument similar to that given above for finite unions).  Since 
$C_{\lambda}$ is ultrafilter closed then we have proven that 
$P_{\calU} \in C_{\lambda}$, for each $\lambda  \in \Lambda$, 
and hence 
$P_{\calU}\in C$.

\end{itemize}

\noindent Note: The above discussion assumes that the sets are 
infinite and the ultrafilters are nonprincipal.  The proofs are 
completely routine otherwise.

\vspace{.2in}

It is natural 
to ask then how this topology compares with the other topologies we 
have defined on $\Spec(R)$.  The goal of this paper is to demonstrate 
that the ultrafilter topology coincides with the patch topology.

It should be noted that related results were obtained recently in  
\cite {CD-N}.

One direction is easy.  We begin with that result.

\begin{prro}  \label{CimpliesU} Let $R$ be a ring and let $C 
\subseteq \Spec(R)$ be a 
collection of prime ideals.  Suppose that $C$ is closed in the patch 
topology.  Then $C$ is also closed in the ultrafilter topology.
\end{prro}

\noindent {\it Proof}. We consider each of the defining classes of 
closed sets for the patch topology separately.

Suppose that $C = V(I)$ for some ideal $I \subseteq R$.  Let $\calU$ 
be a nonprincipal ultrafilter on $C$ and 
construct the prime ideal $P_{\calU}$.  Let $a\in I$.  Then 
$a \in P$ for every prime $P \in C$.  Hence $C \subseteq V(a)$.  
Since our ultrafilter was defined on $C$ we know that $C \in \calU$ 
and so  $V(a) \cap C = C  \in \calU$. 
Hence, $I \subseteq P_{\calU}$.  Since $C$ is defined as the 
collection 
of all primes that contain $I$ this implies that $P_{\calU} \in C$.

Now assume that $C = D(a)$ for some element $a \in R$.   
Let $\calU$ be a nonprincipal ultrafilter on $C$ and construct the 
prime 
ideal $P_{\calU}$.  Then the set $V(a)$  contains none of the primes 
in $C$.  Hence, $\emptyset = V(a) \cap C$ is not 
in $\calU$ since, as noted earlier, no ultrafilter contains the empty set.  
Hence, $a$ is not an element 
of $P_{\calU}$.  It follows 
by definition that $P_{\calU} \in C$.

Finally note that, in both of the above settings, the assertion that $P_{\calU} \in C$ 
is trivial if the 
ultrafilter $\calU$ is principal.
Therefore, we have proven that all of the closed sets which generate 
the patch 
topology are also closed in the ultrafilter topology.  The result 
follows immediately. \hfill $\Box$

\vspace{.2in}

The opposite direction is somewhat harder.  We make use of the 
von Neuman regular ring $T(R)$   canonically  associated with 
the ring 
$R$   (details below).   We start by recalling   a notation 
and   some easy facts concerning the von 
Neumann regular rings.    If $f: R \rightarrow S$ is a ring 
homomorphism, $f^a : \Spec(S) \rightarrow \Spec(R)$ denotes the map 
defined by $f^a(Q) := f^{-1}(Q)$, for each $Q \in \Spec(S)$.

\begin{leem}
The following are equivalent:
 \begin{enumerate}
 \item [\bf (i) \rm]  $R$ is von Neumann regular (i.e. for each $a 
\in R$ there exists $x \in R$ such that $a^2 x =a$).
 
  \item [\bf (ii) \rm]  for each $a \in R$ there exists a unique 
element $a^{(-1)} \in R$ (called \rm the punctual inverse of \it 
$a$)  such that
  $a^2a^{(-1)}  =a$  and $\left(a^{(-1)}\right)^2 a  = a^{(-1)}$.

 \end{enumerate}
 \end{leem}

 \begin{leem}
 \label{lm:2}
 Let $R$ be a von Neumann regular  ring and  set  $e(a) := 
aa^{(-1)}$,   for each $a \in R$.   
  Then the following hold:
 
  \begin{enumerate}
 \item [\bf (1) \rm]  $e(a)$ is idempotent.
 
    
  
  \item [\bf (2) \rm]  $aR = e(a)R$. 
  
    \item [\bf (3) \rm]  If $a, b \in R$, then $
    aR + bR  = e(a)R +e(b)R =\left(e(a) +e(b)(1-e(a)\right)\!R\,.
    $
 \end{enumerate}
 \end{leem}

  \begin{prro}
  \label{von-neumann} Let $R$  be any ring,  let $\{X_a \mid a \in 
R\}$ be a family of indeterminates (one for each element $a \in R$) 
and let $I_R$ be  the ideal generated by $ \{a^2 X_a -a, aX^2_a -X_a 
\mid a \in R\}$ in the polynomial ring  $R[ X_a \mid a \in R]$. Set:
 $$
 T(R) :=\frac{R[ X_a \mid a \in R]}{I_R }\,.
 $$
 Then
  \begin{enumerate}
 \item [\bf (1) \rm]  $T(R)$ is von Neumann regular (thus, in 
particular, every finitely generated ideal of 
$T(R)$ is principal).
 
  \item [\bf (2) \rm]   The canonical embedding $\iota: R \rightarrow 
T(R)$ is an epimorphism.
  
  \item   [\bf (3) \rm] $\iota: R \rightarrow T(R)$ is an isomorphism 
if and only if $R$ is von Neumann regular.
  
    \item   [\bf (4) \rm]   Let \ $\Spec(R)_{\Z}$ \  [respectively, 
\ 
$\Spec(T(R))_{\Z}$] \ be \ $\Spec(R)$ \ [respectively, \  
$\Spec(T(R))$]  \ endowed with the Zariski topology and   let 
$\iota^a : \Spec(T(R)) \rightarrow 
\Spec(R)$ be the canonical map associated to $\iota: R \rightarrow 
T(R)$.     Then   $\iota^a : \Spec(T(R))_{ \Z } \rightarrow 
\Spec(R)_{ \Z }$ is continuous and bijective.

     \item   [\bf (5) \rm] Let \ $\Spec(R)_{\C}$ \  [respectively, \ 
$\Spec(T(R))_{\C}$] \ be \ $\Spec(R)$ \ [respectively, \  
$\Spec(T(R))$]  \ endowed with the patch topology.  Then $\iota^a : 
\Spec(T(R))_{ \Z } \rightarrow \Spec(R)_{\C}$ is an 
homeomorphism. In 
particular,  $\Spec(T(R))_{ \Z } $   coincides with  $ 
\Spec(T(R))_{\C}$. 
    
  \end{enumerate}
  \end{prro}

 Details concerning the 
 von Neumann regular ring $T(R)$    canonically   associated with 
the ring 
$R$  and the proof of the previous (and related) results can be found 
in \cite{Olivier_1}, \cite{Olivier_2}, \cite[Chapter 3, Exercise 
30]{AM}, 
\cite[Chapitre 1, \S 1, Exercices 16, 17, 18; Chapitre 2 \S 4, 
Exercice 16]{B}, 
\cite[ pages 99, Exercice 8 page 119]{Lafon} and \cite{DFP}.

%

%







%



%



%


%



%

%



%



%


%





%



%

For the sake of simplicity, from now on, we identify $R$ with its 
canonical image in $T(R)$.  Consider a subset $C$ of $\Spec(T(R))$. 
Then consider also the 
corresponding collection $C_{R} := \{P \cap R\mid P \in C\}$ of 
primes in $R$.  Note that by Proposition 
\ref{von-neumann} (4),  we can identify the sets 
of primes of $R$ and $T(R)$   by using just the elements that 
lie in $R$.  Define an ultrafilter $\calU$ on $C$.  Then we can define 
a corresponding ultrafilter $\calU_{R}$ on $C_{R}$.  Consider the 
prime $P_{\calU}$ of $T(R)$ defined by $\calU$.  Consider the prime 
$P_{\calU} \cap R$ of $R$.  Choose an element 
$d \in R \subseteq T(R)$.  Define $V_{R}(d)$ to be all the primes of 
$R$ which contain $d$.  Recall that 
$P_{\calU} = \{d \in T(R) \, | \, V(d) \cap C \in \calU \}$.  Then since 
we are assuming that $d \in R$ we have $d \in P_{\calU}$ exactly when  
$V_{R}(d) \cap C_{R} \in \calU_{R}$.  But this is exactly the condition 
necessary for $d$ to lie in $P_{\calU}$.  The point of this is that 
the ultrafilter limit primes of a collection of primes of $T(R)$ will 
correspond precisely to the ultrafilter limit primes of the 
corresponding collection of primes of $R$.  This implies that the 
contraction map $\iota^a : \Spec(T(R)) \rightarrow 
\Spec(R)$ and its inverse send closed sets to closed sets with respect to the ultrafilter topologies.   
The following result is then clear.

\begin{prro} The contraction map $\iota^a : \Spec(T(R)) \rightarrow 
\Spec(R)$ is a homeomorphism  with 
respect to the ultrafilter topologies.  
\end{prro}

So we have homeomorphisms between  $\Spec(T(R))$  and $\Spec(R)$ with 
respect to both the patch topology and the ultrafilter topology.  All 
that remains is to show that the ultrafilter and patch topologies 
coincide on $\Spec(T(R))$.

We can use the result of Proposition \ref{CimpliesU} to show that a 
subset of 
$\Spec(T(R))$ which is closed in the patch topology (or, in the 
Zariski topology    by Proposition \ref{von-neumann} (5)) \ is 
also closed in 
the ultrafilter topology.  We prove the converse.

\begin{prro}  Let $R$ be a ring and let $C$ be a collection of prime 
ideals in $\Spec(T(R))$.  Suppose that $C$ is  closed in the ultrafilter 
topology.  Then $C$ is also closed in the patch topology.
\end{prro}

\noindent {\it Proof}.  We prove the contrapositive. 

Suppose that $C$ is a collection of prime 
ideals in $\Spec(T(R))$which is not closed in the patch  
topology. 
 Let $I := \cap_{P \in C}P \ (\subseteq 
T(R))$.  Since $C$ is not 
closed in the patch topology (or, equivalently, in the Zariski 
topology)  of $\Spec(T(R))$ then 
$V(I)$ is properly larger than $C$.  
Let $ \overline{P}$ be a prime ideal in $T(R)$ such that $I 
\subseteq  \overline{P}$ but 
$ \overline{P}  \notin C$.   Choose a   nonzero   element 
$a \in  \overline{P}$.  We claim that 
$a \in Q$,  for some $Q \in C$.  Suppose not.  We can find an 
element $x \in T(R)$ such that $a^{2}x - a = 0$.  Note then that 
$a(ax - 1) = 0$.  It follows that any prime which does not contain 
$a$ must contain $ax - 1$.  Hence, $ax- 1$ is contained in every 
prime in $C$.  This implies that $ax- 1$ is contained in $I$.  This 
then implies that $ax- 1$ is contained in $ \overline{P}$.  
Since $a$ and 
$ax - 1$ are relatively prime and are both contained in 
${ \overline{P}}$ we have a 
contradiction.  

For each element $a \in  \overline{P}$ let $V_{C}(a)$ represent 
all of the primes of $C$ which contain $a$ (i.e., $V_C(a) = V(a) \cap 
C$).  (Note that we have proven in the preceding paragraph that 
$V_{C}(a)$ is not empty.)  Our goal is to build an 
ultrafilter on $C$ such that each set $V_{C}(a)$ is in the 
ultrafilter.  If we can show that a collection of nonempty subsets of a set is 
closed under finite intersection then we will have 
shown that all subsets in the collection lie in a filter.  Then 
note that the collection of all filters on a set can be partially 
ordered under inclusion.  Zorn's Lemma implies that maximal filters 
exist and that every filter is contained in a maximal filter.  The 
maximal filters correspond precisely to the ultrafilters. 
So our goal is to prove that the collection of 
sets of the form $V_{C}(a)$, with $a \in \overline{P}$, is 
closed under finite 
  intersections.   
This is equivalent to showing that if $J$ is a  nonzero   
finitely generated 
ideal contained in $\overline{P}$ then there is some prime $Q 
\in C$ such 
that $J \subseteq Q$.  Recall however, that in $T(R)$ all finitely 
generated ideals are principal (Proposition \ref{von-neumann} (1)).  
It follows that such an 
ideal $J$ is actually principal and 
since the generator must lie in $ \overline{P}$ we have already 
proven what we 
need.  (Note that we have proven that a finite intersection of ideals of the form 
$V_{C}(a)$ has that same form and hence is not empty.) 
So we let $\calU   :=  \calU{(\overline{P})}$ be an 
ultrafilter which contains 
all of the 
sets $V_{C}(a)$ for $a \in  \overline{P}$.  Then we construct 
the ultrafilter 
limit prime $P_{\calU}$.  The construction was designed so that 
$ \overline{P} \subseteq P_{\calU}$.  However, $T(R)$ is a 
zero-dimensional ring 
so 
we have actually proven that $ \overline{P} = P_{\calU}$.  Since 
we assumed that 
$ \overline{P}$ 
is not in $C$ this implies that $C$ is not ultrafilter closed. \hfill 
$\Box$

\vspace{.2in}

The preceding result finishes the last step in our main theorem.

\begin{thee} Let $R$ be a ring.  Then the patch topology and the 
ultrafilter topology on the collection $\Spec(R)$ of prime ideals of 
$R$ are identical.
\end{thee}

\vspace{8pt}

\noindent \bf Added in proof.  \rm  The authors have become aware recently that Gabriel Picavet obtained some results regarding ultrafilters and spectral topologies which appeared in:
{\it Ultrafiltres sur un espace spectral, anneaux de Baer, anneaux a spectre minimal compact}, Math. Scand, {\bf 46}, (1980), 23-53.  The authors also wish to thank Picavet for pointing  out that the results of J.-P. Olivier cited in this paper can also be found in the paper by Olivier published as {\it L'anneau absolument plat universel, les epimorphismes et les parties constructibles} Boletin de la Sociedad Matematica Mexicana, {\bf 23} (1978), 68-74.

\end{document}